\documentclass[10pt,a4paper,leqno]{amsart}

\usepackage{amsmath,amssymb,amsfonts,amsthm,epsfig}

\def\dis
{\displaystyle}

\def\R{{\mathbb R}}
\def\N{{\mathbb N}}

\def\S{{\mathcal S}}

\def\virgp{\raise 2pt\hbox{,}}

\def\({\left(}
\def\){\right)}
\def\<{\left\langle}
\def\>{\right\rangle}

\def\Eq#1#2{\mathop{\sim}\limits_{#1\rightarrow#2}}
\def\Tend#1#2{\mathop{\longrightarrow}\limits_{#1\rightarrow#2}}

\def\d{{\partial}}
\def\e{\varepsilon}
\def\si{{\sigma}}
\def\O{{\mathcal O}}
\def\om{{\omega}}

\def\w{{\tt w}}
\def\1{{\rm 1\mskip-4.5mu l} }

\theoremstyle{plain}
\newtheorem{theo}{Theorem}[section]
\newtheorem{lem}[theo]{Lemma}

\newtheorem{prop}[theo]{Proposition}

\theoremstyle{definition}

\theoremstyle{remark}
\newtheorem*{rema}{Remark}

\numberwithin{equation}{section}

\begin{document}

\title[NLS with potential]{Linear vs. nonlinear effects for nonlinear
Schr\"odinger 
equations with potential}  
\author[R. Carles]{R{\'e}mi Carles}
\email{remi.carles@math.univ-rennes1.fr}
\thanks{The author acknowledges support by the European network HYKE,
funded  by the EC as contract  HPRN-CT-2002-00282}
\begin{abstract}
We review some recent results on nonlinear Schr\"odinger
equations with potential, with emphasis on the case where the
potential is a second order polynomial, for which the interaction
between the linear dynamics caused by the potential,
and the nonlinear effects, can be described quite
precisely. This includes semi-classical
r\'egimes, as well as finite time blow-up  and scattering issues. We
present the 
tools used for these problems, as well as their limitations, and
outline the arguments of the proofs.
\end{abstract}

\maketitle

\section{Introduction}
\label{sec:intro}

This paper is a survey of some recent results on nonlinear
Schr\"odinger equations with potential. A particular attention is paid
to the case where the potential is a second order polynomial. In this
case, the fundamental solution is known explicitly, through a
generalized Mehler's formula; the linear dynamics is well
understood. The most important remark is that we can assess the action
of some Heisenberg
observables (which can be exactly computed) on a class of
nonlinearities. With this tool, we can understand the interaction
between the linear effects caused by the potential, and the nonlinear
effects. In the semi-classical r\'egime $\e \to 0$, we emphasize
some critical 
scales measuring this interaction, and describe the critical
phenomena. When $\e=1$, we observe the
effect of such potentials on  finite 
time blow-up; the potential may ``create'' some blow-up (think of the
harmonic potential), or delay, even prevent, this phenomenon
(``repulsive'' harmonic potential). 

We give two motivations to study nonlinear
Schr\"odinger equations with potential. The first one arises from
physics, for Bose--Einstein condensates. The model
involves the equation
\begin{equation}\label{eq:BEC}
i\hbar \d_t u^\hbar +\frac{1}{2}\hbar^2 \Delta u^\hbar = V(x)u^\hbar +
a\hbar^2 |u^\hbar|^{2\si}u^\hbar\ ;\quad x\in \R^n\, ,\ a\in \R\, ,\
\si \in \N\, , 
\end{equation}
where the role of the potential $V$ is to confine particles. The 
cases most currently considered are when $V$ is quadratic (isotropic
or anisotropic harmonic potential),  when $V$ is lattice periodic, or
when $V$ is the sum of two such potentials (see
e.g. \cite{PiSt,BSTH,KNSQ,DFK}). We shall not discuss the physical
relevance of this model, but notice that explicit mathematical
formulae are available in the 
case of the harmonic potential. Our main motivation is more from a
mathematical point of 
view. The presence of a potential alters the propagation of the wave
in the linear case. Similarly, nonlinear problems may lead to typical
phenomena, such as finite time blow-up. How can these two effects
interact? This vague question seems to raise many complex issues. We
give some very partial answers, essentially restricted to the case
where the potential is a second order polynomial. These examples may
be viewed as a first step for a general study, supporting or
contradicting the intuition. This may be compared with the approach of
N.~Burq, P.~G\'erard and N.~Tzvetkov, who analyze the role of geometry
for nonlinear Schr\"odinger equations on compact manifolds (see
e.g. \cite{BGT,BGTMRL}). 

The initial value problem we study is
\begin{equation}\label{eq:nlsp}
i\e\d_t u^\e +\frac{1}{2}\e^2\Delta u^\e = V(x)u^\e+ \lambda^\e
|u^\e|^{2\si}u^\e\quad ; \quad 
u^\e\big|_{t=0} = u_0^\e\, ,
\end{equation}
where $\e\in]0,1]$, $x\in \R^n$, $V\in C^\infty(\R^n;\R)$,
$\lambda^\e\in\R$, and $\si >0$ with $\si<\frac{2}{n-2}$ if $n\ge
3$ (the nonlinearity is $H^1$ sub-critical). We consider two r\'egimes
for the parameter $\e$: 
\begin{itemize}
\item The semi-classical limit $\e\to 0$. This r\'egime gives hints
to understand high-frequency phenomena and provides us with tools which
are extremely useful in the nonlinear case, when $V$ is a second order
polynomial.
\item The case $\e=1$. This case is better
understood thanks to the semi-classical analysis. We study in
particular the global existence issue and two of its companions: finite time
blow-up and scattering. 
\end{itemize}
We shall not discuss any regularity issue here, and always assume
$u_0^\e \in \Sigma$, where
\begin{equation*}
\Sigma := \left\{ f \in \S'(\R^n)\ ;\ \left\| f\right\|_\Sigma
:= \|f\|_{L^2(\R^n)} + \|\nabla f\|_{L^2(\R^n)}+ \|x
f\|_{L^2(\R^n)}<+\infty \right\}\, . 
\end{equation*}
We do not discuss the question of solitons either. 
We denote 
\begin{equation*}
H_V^\e =-\frac{1}{2}\e^2 \Delta +V(x)\  ; \  U_V^\e(t) =
e^{-i\frac{t}{\e}H_V^\e}\, ,\ \text{and simply $H_V$ and $U_V$ when $\e=1$.}
\end{equation*}
This paper is organized as follows. Section~\ref{sec:general} is
devoted to general results. First, we recall some classical results in
the case $V\equiv 0$ and $\e=1$, as well as some techniques to prove
them. We then discuss which ones can be easily generalized when $V$ is
not identically zero, and present some cases where a change of
variables makes it possible to relate the case $V\equiv 0$ to the
case where $V$ is not trivial. 

In Section~\ref{sec:mehler}, we recall
the generalized Mehler's formula, and motivate the introduction of
some particular Heisenberg observables, from a nonlinear point of
view. 

In Section~\ref{sec:semiclas}, we state some results about
\eqref{eq:nlsp} in the semi-classical limit; we outline the
techniques, and discuss their limitations. 

In
Sections~\ref{sec:blowup} and \ref{sec:global}, we set $\e=1$, and
analyze the influence of the potential $V$ on finite time blow-up and
global existence issues. 
\section{General setting and consequences}
\label{sec:general}

\subsection{Some results on the nonlinear Schr\"odinger
equation}\label{sec:nls}
In this paragraph, we assume $\e=1$, and recall a few things about the
Cauchy problem
\begin{equation}\label{eq:nls}
i\d_t u +\frac{1}{2}\Delta u =\lambda
|u|^{2\si}u\quad ; \quad 
u\big|_{t=0} = u_0\, ,
\end{equation}
where $x\in \R^n$, 
$\lambda\in\R$, and $\si >0$ with $\si<\frac{2}{n-2}$ if $n\ge 3$. All
the results we mention can be found in \cite{CazCourant}. 
\begin{theo}\label{theo:nlslocal}
Suppose $u_0\in H^1(\R^n)$, and $\si>0$ as above. Then there exist
$T_*,T^*>0$,
and a unique maximal solution 
$$u\in C\(]  -T_*,T^*[;H^1\)\cap L^q_{\rm loc}\(]-T_*,T^*[;
W^{1,2\si+2}\)$$
to \eqref{eq:nls}, where $q=\frac{4\si+4}{n\si}$. It is maximal in the
sense that if $T^*<\infty$, then
\begin{equation*}
\lim_{t\to T^*}\|\nabla_x u(t)\|_{L^2}=+\infty\, .
\end{equation*}
In addition, the following quantities are independent of time:
\begin{equation}\label{eq:conservations}
\begin{aligned}
&\text{Mass: }\left\| u(t)\right\|_{L^2} \equiv \left\|
u_0\right\|_{L^2}\, ,\\
&\text{Energy: } E_0 := \frac{1}{2}\left\| \nabla_x
u(t)\right\|_{L^2}^2 +\frac{\lambda}{\si +1}\left\|
u(t)\right\|_{L^{2\si +2}}^{2\si +2}
=\text{const.}
\end{aligned}
\end{equation}
\end{theo}
This result was first proved in \cite{GV79Cauchy} (with a slightly
different statement), and revisited in \cite{GV85c,Yajima87}. The
proof relies on a fixed point argument on Duhamel's formula
\begin{equation}\label{eq:duhamel}
u(t)=U_0(t)u_0 -i\lambda \int_0^t U_0(t-s)\(|u|^{2\si}u\)(s) ds\, ,
\end{equation}
where following the notation introduced in
Section~\ref{sec:intro}, $U_0(t) = e^{i\frac{t}{2}\Delta}$. The modern
tool to prove this result is \emph{Strichartz estimates}, after
\cite{Strichartz}. We recall a statement for the case of Schr\"odinger
equations \cite{Yajima87,KT} (see \cite{CaQuad} for the following adaptation):
\begin{lem}\label{lem:strichartz}
Let $(U(t))_{t\in\R}$ be a unitary group on $L^2(\R^n)$, satisfying
the dispersive estimate 
\begin{equation}\label{eq:dispgen}
\left\|U(t)\right\|_{L^1\to 
L^\infty}  
\le \w(t)^{n/2}\, ,\quad \text{with }\w\ge 0\ \text{and }\w
\in L^1_w(\R)\, .
\end{equation}
We say that a pair $(q,r)$ is \emph{admissible} if
\begin{equation*}
\frac{2}{q}+\frac{n}{r}=\frac{n}{2}\, ;\quad  q,r\ge 2,\ \text{ with
}r<\frac{2n}{n-2}\, \cdot
\end{equation*} 
Then for any $T\in
\overline{\R}_+$, and any 
admissible pairs $(q,r)$ and $(\widetilde q,\widetilde r)$,
we have
\begin{align*}
\left\| U(t)f\right\|_{L^q(]-T,T[;L^r)} &\le
\left\|\w\1_{]-2T,2T[}\right\|_{L^1_w}^{1/q} \left\| 
f\right\|_{L^2}\, ,\\
\left\| \int_0^t U(t-s)F(s)ds \right\|_{L^q(]-T,T[;L^r)} &\le  
\left\|\w\1_{]-2T,2T[}\right\|_{L^1_w}^{1/q +1/\widetilde q}
\left\| F\right\|_{L^{\widetilde q'}(]-T,T[;L^{\widetilde
r'})}\, .
\end{align*}
\end{lem}
It is straightforward that $U_0(t)$ satisfies the assumptions of
Lemma~\ref{lem:strichartz}, with $\w(t) = (2\pi |t|)^{-1}$, from the
formula
\begin{equation}\label{eq:U0}
U_0(t)f (x) = \frac{1}{(2i\pi t)^{n/2}} \int_{\R^n}
e^{i\frac{|x-y|^2}{2t}}f(y)dy\, .
\end{equation}
Theorem~\ref{theo:nlslocal} then follows from Strichartz and
H\"older's inequalities, as well as Gagliardo--Nirenberg inequalities
when $\si\ge \frac{2}{n}$. 
A similar result holds when $u_0\in \Sigma$; essentially, we can
remember that in addition, $u\in C(\R;\Sigma)$.
 
In the $H^1$ or $\Sigma$ case, global existence  can be deduced from
the conservations 
\eqref{eq:conservations}, when $\lambda >0$ (defocusing nonlinearity),
or $\lambda <0$ and $\si<\frac{2}{n}$ for instance. If $\lambda<0$ and
$\si\ge \frac{2}{n}$, finite time blow-up may occur:
\begin{prop}[Virial Theorem, \cite{Glassey}]\label{prop:glassey}
Assume $\lambda<0$ and
$\si\ge \frac{2}{n}$ (with $\si<\frac{2}{n-2}$ if $n\ge 3$). If
$u_0\in \Sigma$ is such that $E_0<0$, then solutions to \eqref{eq:nls}
blow up in finite time, in the future and in the past: $T_*$ and $T^*$
are finite. 
\end{prop}
The idea of the proof for this result is to introduce the function
$y(t)=\|xu(t)\|_{L^2}^2$ and to notice that under the above
assumptions, $\ddot y(t)\lesssim E_0<0$. Since $y(t)\ge 0$, this
proves that a singularity appears, both for $t>0$ and for $t<0$. 

A slightly different approach to recover this result is to use the
Galilean operator $J(t)=x+it\nabla_x$, and the
\emph{pseudo-conformal conservation law}, discovered in
\cite{GV79Scatt}:
\begin{equation}\label{eq:pseudoconf}
\frac{d}{dt}\left(\frac{1}{2}\| J(t)u
\|^2_{L^2}+\frac{\lambda}{\si+1}t^2\|u
    (t)\|^{2\si+2}_{L^{2\si +2}}\right)=\frac{\lambda}{\si+1}(2-n\si)t\|u
    (t)\|^{2\si+2}_{L^{2\si+2}}.
\end{equation}
As noticed in \cite{Weinstein83}, expanding the above formula, one
retrieves Proposition~\ref{prop:glassey}. Another application of the
pseudo-conformal conservation law, as motivated in \cite{GV79Scatt},
concerns scattering theory:
\begin{prop}[Scattering theory in $\Sigma$,
\cite{GV79Scatt,HT87,CW92,NakanishiOzawa}]\label{prop:scattnls}  
Assume 
$\frac{2}{n+2}<\sigma<\frac{2}{n-2}$ if $n\geq 2$, $\sigma >1$ if
$n=1$. Assume in addition $\lambda>0$. Denote 
\begin{equation*}
\sigma_0(n):=\frac{2-n +\sqrt{n^2+12n +4}}{4n}\, \cdot
\end{equation*}
\begin{itemize}
\item Let $u_-\in \Sigma$. There exists a unique $u \in 
C(\R;\Sigma)$ solution of \eqref{eq:nls}, such that 
\begin{equation*}
\lim_{t\rightarrow -\infty}\|U_0(-t)u(t)-u_-\|_\Sigma
=0\, .
\end{equation*}
\item Let $u_0\in \Sigma$. If
$\sigma\ge \sigma_0(n)$ or if $\|u_0\|_\Sigma$ is sufficiently
small, then  there exists a unique $u_+ \in \Sigma$ such that if  $u \in 
C(\R;\Sigma)$ is the solution of \eqref{eq:nls}, then
\begin{equation*}
\lim_{t\rightarrow +\infty}\|U_0(-t)u(t)-u_+\|_\Sigma =0\, .
\end{equation*}
\end{itemize}
We denote by $S:u_-\mapsto u_+$ the scattering operator. 
\end{prop}
We emphasize two properties of the Galilean operators, from which
the above results (Equation~\eqref{eq:pseudoconf} and
Proposition~\ref{prop:scattnls}) follow:
\begin{itemize}
\item[(i)] $J(t)$ is an Heisenberg observable: 
\begin{equation}\label{eq:JHeisenberg}
J(t) =U_0(t)xU_0(-t)\, .
\end{equation}
\item[(ii)] It also reads 
\begin{equation}\label{eq:Jfactor}
J(t)= i t \, e^{i\frac{|x|^2}{2t}}\nabla_x \( e^{-i\frac{|x|^2}{2t}}\,
\cdot\)\, .
\end{equation}
\end{itemize}
As a consequence of \eqref{eq:JHeisenberg}, $J(t)$ commutes with the
linear part of \eqref{eq:nls},
\begin{equation}\label{eq:Jcommute}
\left[ J(t),i\d_t +\frac{1}{2}\Delta\right]=0\, .
\end{equation}
The second point has two interesting straightforward consequences:
\begin{itemize}
\item[(ii)'] Weighted Gagliardo--Nirenberg inequalities: for
$2\le p <\frac{2n}{n-2}$ ($2\le p\le \infty$ if $n=1$), there exists
$C_p$ depending only on $n$ and $p$ such that 
\begin{equation}\label{eq:GNlibre}
\left\| f \right\|_{L^p(\R^n)}\le \frac{C_p}{|t|^{\delta(p)}} \left\|
f \right\|_{L^2(\R^n)}^{1-\delta(p)} \left\|
J(t) f \right\|_{L^2(\R^n)}^{\delta(p)}\ ; \quad \delta(p):=n\(
\frac{1}{2}-\frac{1}{p}\)\, .
\end{equation}
\item[(ii)''] If $F(z)=G(|z|^2)z$ is $C^1$, then $J(t)$
acts like a derivative on $F(w)$:
\begin{equation}\label{eq:Jder}
J(t)\(F(w)\) = \d_z F(w)J(t)w -\d_{\overline z} F(w)\overline{ J(t)w
}\, . 
\end{equation}
\end{itemize}
Roughly speaking, \eqref{eq:Jcommute} and \eqref{eq:Jder} make it
possible to have estimates for $J(t)u$ in the same way as for
$\nabla_xu$, when $u$ solves \eqref{eq:nls}. Then \eqref{eq:GNlibre}
yields dispersive estimates for the nonlinear equation, which are the
same as in the linear case, provided that $J(t)u\in
L^\infty_tL^2_x$. These arguments are the key ingredients for the
first point of Proposition~\ref{prop:scattnls} (existence of wave
operators). Then \eqref{eq:pseudoconf} yields estimates on
$\|J(t)u\|_{L^2}$ which, along with the conservation of mass
\eqref{eq:conservations}, prove the second point of
Proposition~\ref{prop:scattnls} (asymptotic completeness). 
\begin{rema}
The properties \eqref{eq:Jcommute}, \eqref{eq:GNlibre} and
\eqref{eq:Jder} are the analog of those satisfied by the conformal
Killing vector-fields used to study the 
wave equation (\cite{Klainerman85}). 
\end{rema}

\subsection{Introducing a potential}

We now turn to the case of \eqref{eq:nlsp}. Drawing a parallel with
the case of \eqref{eq:nls}, the first tool we seek is Strichartz
estimates for $U_V$. This has been, and this is still,  a very active
 area of research. The possible presence of eigenvalues shows that in
general, one cannot expect not only the same global dispersion as for
$U_0$, but also global in time Strichartz estimates. In a more subtle
way, resonances may also be an obstacle to dispersion
(\cite{JensenKato}). Many results have
been obtained though, and we refer to the introduction of
\cite{RodnianskiSchlag} for a very nice review. In
\cite{RodnianskiSchlag}, Strichartz estimates are obtained for
potentials which may depend of time, a case we do not consider. Notice
that these involved results rely on perturbation arguments (see
\cite{NierSoffer} for an
interesting exception). We would
like to consider the case where the potential may really change the
dynamics of the Laplacian, and as introduced in \eqref{eq:nlsp}, we
assume that $V\in C^\infty(\R^n;\R)$. 

It is well known that if $V(x)\ge -a |x|^2-b $ for some $a,b>0$, then
$H_V$ is essentially self-adjoint on $C_0^\infty(\R^n)$, and that this
is a sharp result: if $V(x)=-|x|^4$ for instance, this property
fails (classical trajectories can reach infinite speed, see
\cite{ReedSimon2,Dunford}). On the other hand, if $V$ is positive and
grow faster than quadratically, then at least for $n=1$, 
the kernel of $U_V$ is
nowhere $C^1$ (\cite{Yajima96}), but smoothing properties make it
possible to solve the nonlinear problem \eqref{eq:nlsp} in some
cases (\cite{YajZha01}). We shall restrict our attention to
\emph{sub-quadratic potentials}, for which 
it is possible to construct a parametrix, locally in time:
\begin{prop}[\cite{Fujiwara79,Fujiwara}]\label{prop:parametrice}
Let  $V\in C^\infty(\R^n;\R)$ be such that $\d^\alpha V\in
L^\infty(\R^n)$ for $|\alpha|\ge 2$. Then there exists
$\delta >0$ independent of $\e\in ]0,1]$ 
such that for $|t|\le \delta$, 
\begin{equation}\label{eq:solfond}
U_V ^\e(t)f(x)= e^{-in\frac{\pi}{4}\operatorname{sgn}t}
\frac{1}{\left|2\pi \e t \right|^{n/2}}\int_{\R^n}
k^\e(t,x,y) e^{iS(t,x,y)/\e}  f(y)dy \, ,
\end{equation}
where $S$ solves the eikonal equation
\begin{equation}\label{eq:eikonale}
\d_t S +\frac{1}{2}|\nabla_x S|^2 +V(x)=0\, ,
\end{equation}
and  $\d_x^\alpha \d_y^\beta k^\e\in L^\infty( ]-\delta,\delta[\times
 \R^{2n})$ for all $\alpha, \beta\in \N^n$, 
uniformly in $\e\in ]0,1]$.
\end{prop}
This result relies on perturbation arguments: for small
times, the influence of the potential on classical trajectories is
controlled, hence a formula similar to \eqref{eq:U0}. 

Since we know that $U_V^\e$ is a unitary group on $L^2$, the above result
shows that for $|t|\le \delta$, $U_V^\e$ is dispersive:
\begin{equation*}
\left\| U_V^\e(t)\right\|_{L^1\to L^\infty}\lesssim
\frac{1}{|\e t|^{-n/2}}\, \cdot
\end{equation*}
We infer that there are local in time Strichartz estimates. Notice
that in general, one cannot expect global in time estimates: when $V$
is the harmonic potential or is lattice periodic, $H_V$ (hence
$H_V^\e$) has eigenvalues. 

One can then mimic the proof of Theorem~\ref{theo:nlslocal}. Unlike
for the case of \eqref{eq:nls}, the gradient does not commute with the
linear part of the equation:
\begin{equation*}
\(i\e\d_t +\frac{1}{2}\e^2\Delta -V(x)\)\nabla_x u^\e= u^\e\nabla_x V
+ \lambda^\e \nabla_x\( |u^\e|^{2\si}u^\e\)\, .
\end{equation*}
The new term is $u^\e\nabla_x V$. Since $V$ is sub-quadratic, $|\nabla
V(x)|\lesssim 1 +|x|$, which suggests to consider $xu^\e$ as a third
unknown, after $u^\e$ and $\nabla_x u^\e$. We then have
\begin{equation*}
\(i\e\d_t +\frac{1}{2}\e^2\Delta -V(x)\)x u^\e= \nabla_x u^\e
+ \lambda^\e x\( |u^\e|^{2\si}u^\e\)\, .
\end{equation*}
It is easy to prove:
\begin{lem}\label{lem:exlocV}
Suppose that $V\in C^\infty(\R^n;\R)$ is sub-quadratic, $u_0^\e\in
\Sigma$, $\lambda^\e\in \R$, and $\si>0$ with $\si<\frac{2}{n-2}$ if
$n\ge 3$. Then there exist $T^\e>0$
(depending on $n$, $\si$, $V$, $|\lambda^\e|$ and $\|u_0^\e\|_\Sigma$),
and a unique solution 
$$u^\e\in C\(]  -T^\e,T^\e[;\Sigma \)\cap L^q_{\rm loc}\(]-T^\e,T^\e[;
W^{1,2\si+2}\)$$
to \eqref{eq:nlsp}, where $q=\frac{4\si+4}{n\si}$. 
The following quantities are independent of time:
\begin{equation*}
\begin{aligned}
&\text{Mass: }\left\| u^\e(t)\right\|_{L^2} \equiv \left\|
u_0^\e\right\|_{L^2}\, ,\\
&\text{Energy: } E_V^\e := \frac{1}{2}\left\| \e\nabla_x
u^\e(t)\right\|_{L^2}^2 +\frac{\lambda^\e}{\si +1}\left\|
u^\e(t)\right\|_{L^{2\si +2}}^{2\si
+2}+\int_{\R^n}V(x)|u^\e(t,x)|^2dx\, .
\end{aligned}
\end{equation*}
In addition, if $\si<\frac{2}{n}$ ($L^2$ sub-critical case), then one
can take $T^\e=+\infty$. 
\end{lem}
Notice that without further assumption on $V$, we cannot state a
criterion for the obstruction to global existence, as in
Theorem~\ref{theo:nlslocal}. Global existence in the $L^2$
sub-critical case follows from the same ideas as in \cite{TsutsumiL2}:
one has global existence at the $L^2$ level, from Strichartz
inequalities (local in time inequalities are sufficient) and the
conservation of mass. One deduces global 
existence in $\Sigma$ by considering $|u^\e|^{2\si}$ as a potential in
the equations for $\nabla_x u^\e$ and $x u^\e$. 

If in addition $V$ is non-negative, then the assumption $u_0^\e\in
\Sigma$ can be weakened: it suffices to consider initial data in
the domain of $\sqrt{-\Delta +V}$ (see \cite{Oh,CazCourant}). In that
case, the obstruction to global existence is
the same as in Theorem~\ref{theo:nlslocal}. 

\subsection{Removing the potential}
It turns out that for some specific potentials, an explicit change of
variables makes it possible to relate the solutions of \eqref{eq:nlsp}
to the solutions of the same equation with $V\equiv 0$. These
potentials are:
\begin{itemize}
\item The linear potential (Stark potential): $V(x)=E\cdot x$.
\item The isotropic harmonic potential $V(x)=|x|^2$, in the case where
the nonlinearity is $L^2$-critical, $\si=\frac{2}{n}$. 
\item The isotropic \emph{repulsive} harmonic potential $V(x)=-|x|^2$,
in the case where the nonlinearity is $L^2$-critical.
\end{itemize}
Introduce $v^\e$ the solution to the equation with no potential:
\begin{equation*}
i\e\d_t v^\e +\frac{1}{2}\e^2\Delta v^\e = \lambda^\e
|v^\e|^{2\si}v^\e\quad ; \quad 
v^\e\big|_{t=0} =u_0^\e \ \(= u^\e\big|_{t=0} \)\, .
\end{equation*}
\subsubsection{Linear potential}\label{sec:stark}
Assume that $V(x)=E\cdot x$, where $E\in \R^n$ is constant. Then as
noticed in \cite{CN}, the Avron--Herbst formula, discovered in the
linear case $\lambda^\e=0$ (\cite{AvronHerbst}), relates $u^\e$ and
$v^\e$: 
\begin{equation}\label{eq:chgtstark}
u^\e(t,x)= v^\e\left(t, x +\frac{t^2}{2}E \right)e^{-i\left( 
tE\cdot x +\frac{t^3}{6}|E|^2\right)/\e}.
\end{equation}
Therefore, one can use Theorem~\ref{theo:nlslocal} to deduce local
existence results in $H^1$ and see that the only obstruction to global
existence is the unboundedness of $\nabla_x u^\e$ in $L^2$. Similarly,
the Stark potential does not change the possible blow-up time, nor the
scattering 
theory (replace $U_0$ by $U_V$ in Proposition~\ref{prop:scattnls}). In
the case of finite time blow-up, the linear potential only shifts the
set where this phenomenon occurs: if $v^\e$ blows up on a set $X$ at time 
$T^\e>0$ (see e.g. \cite{MerleCMP} for a definition of such a set),
then $u^\e$ blows up (at time $T^\e$) on $X -\frac{(T^\e)^2 }{2}E$.
\subsubsection{Isotropic harmonic potential and conformal
nonlinearity}\label{sec:isoharmocrit}
Assume that $\si=\frac{2}{n}$ and
$V(x)=\frac{\om^2}{2}|x|^2$ with
$\om>0$. Then provided that the right hand side is defined,
\begin{equation}\label{eq:chgtharmo}
u^\e(t,x)= \frac{1}{(\cos \omega t)^{n/2}}e^{-i\frac{\omega}{2\e}x^2
\tan \omega 
t} v^\e \left(\frac{\tan \omega t}{\omega}\virgp\frac{x}{\cos \omega t}
\right). 
\end{equation}
This was first noticed in \cite{Niederer} for the linear case
($\lambda^\e =0$), and in 
\cite{Rybin,CaM3AS} for the nonlinear case with critical
nonlinearity. Note that despite this expression, $u^\e$ needs not
blow-up at time $t_1=\frac{\pi}{2\om}$ (when the cosine is
zero). Indeed, as $t\to t_1$, the time variable $ \frac{\tan \omega
t}{\omega}$ for $v^\e$ goes to infinity. Then dispersive properties of
solutions to \eqref{eq:nls} may compensate the cancellation of the
cosine. 

However, it is true that the harmonic potential generates more
blowing-up solutions (the criterion for finite time blow-up is the
same as in Theorem~\ref{theo:nlslocal}). First, if $v^\e$ blows up at
time $T^\e>0$, then $u^\e$ blows up at time $\frac{\arctan(\om
T^\e)}{T^\e}$, which is always smaller than $T^\e$. Second, if $v^\e$
is not a dispersive solution, typically of solitary wave $v^\e(t,x) =
e^{i\kappa t/\e} Q(\frac{x}{\e})$, then $u^\e$ blows up at time
$\frac{\pi}{2\om}$. From a heuristic point of view, the confining
properties of the harmonic potential are sufficient to concentrate an
energy which is not naturally dispersed. 

Similarly, because the harmonic potential prevents the solutions from
being dispersive as in the case $V\equiv 0$, no scattering theory must
be expected. 

\subsubsection{Isotropic repulsive harmonic potential and conformal
nonlinearity}\label{sec:isorepharmocrit}
Assume that $\si=\frac{2}{n}$ and
$V(x)=-\frac{\om^2}{2}|x|^2$ with
$\om>0$. Then as noticed in \cite{CaSIMA}, 
\begin{equation}\label{eq:chgtrep}
u^\e(t,x)= \frac{1}{\left(\cosh (\omega t)
\right)^{n/2}}e^{i\frac{\omega}{2\e}|x|^2 \tanh (\omega
t)} v^\e\left(\frac{\tanh (\omega t)}{\omega}\virgp \frac{x}{\cosh (\omega t)}
\right).
\end{equation}
A criterion for global existence is not obvious in this case, since
the potential is unbounded from below. It was proved in \cite{CaSIMA}
that for an 
isotropic repulsive harmonic potential (and a general nonlinearity as
in \eqref{eq:nlsp}), the obstruction to global existence is the same
as in Theorem~\ref{theo:nlslocal} (see also
Section~\ref{sec:blowup}). Opposite to the 
harmonic potential mentioned above, the repulsive harmonic potential
tends to prevent blow-up. If $v^\e$ blows up at time $T^\e>0$, then
$u^\e$ blows up at time $\frac{\arg \tanh (\om T^\e)}{\om}$ \emph{if}
$\om T^\e <1$. That means that for ``small'' values of $\om$, blow-up
is delayed. On the other hand, if $\om T^\e \ge 1$, then $u^\e$ does
not ``see'' the blow up of $v^\e$, and is global. If $\om T^\e > 1$,
$u^\e$ has even exponential decay as $t\to +\infty$. The limiting case
$\om T^\e =1$ is interesting: it is known that there exists no ($L^2$)
solitary wave (\cite{JohnsonPan,KavianWeissler,MerleRaphaelGAFA}). The
analog of a solitary wave is a solution that grows exponentially in
time (see \cite{CaSIMA}). 

Unlike the harmonic potential, the repulsive harmonic potential tends
to spread out a solution to \eqref{eq:nls}: we saw how it changes the
finite time blow-up phenomenon. It is not surprising that things go
very well as far as scattering is concerned, thanks to exponential
decay. We study this issue more precisely in Section~\ref{sec:global}. 
\section{Mehler's formula and applications}
\label{sec:mehler}
When $V$ is a
second order polynomial, a lot can be said. First,
the linear case is extremely favorable, because $U_V^\e$ is known
explicitly: this is the (generalized) \emph{Mehler's formula}
(\cite{Feyn,HormanderQuad}). If $V(x) =\sum a_{jk}x_j x_k + \sum b_j x_j
+c$, then reducing the quadratic part, a change of orthonormal basis (which
leaves the Laplace operator invariant) and a change of the origin
simplify the expression of 
$V$:
\begin{equation*}
V(x) = \sum_{j=1}^n \(\delta_j \frac{\om_j^2}{2}x_j^2 +\widetilde b_j
x_j \) +\widetilde c\quad ;\quad \delta_j\in \{-1,0,1\},\ \om_j >0,\
\delta_j\widetilde b_j =0\, . 
\end{equation*}
Using Avron--Herbst formula \eqref{eq:chgtstark}, we can get rid of
linear terms $\widetilde b_j x_j$. Taking $e^{i\widetilde c t/\e}u^\e$
as a new unknown function also removes the constant term. We therefore
assume that $V$ is of the  form
\begin{equation}\label{eq:V}
V(x) = \sum_{j=1}^n \delta_j \frac{\om_j^2}{2}x_j^2 \quad ;\quad
\delta_j\in \{-1,0,1\},\ \om_j >0\, . 
\end{equation}
The Mehler's formula then reads:
\begin{equation}\label{eq:mehler}
U_V^\e(t)f(x)=
\prod_{j=1}^n
\left(\frac{1}{2i\pi\e g_j(t)}\right)^{1/2}
\int_{\R^n}e^{iS(t,x,y)/\e }f(y)dy\, ,
\end{equation}
where
\begin{equation*}
S(t,x,y)= \sum_{j=1}^n \frac{1}{g_j(t)}\(
\frac{x_j^2 +y_j^2}{2}h_j(t) -x_j
y_j\),
\end{equation*}
and the functions $g_j$ and $h_j$, related to the classical
trajectories, are given by:
\begin{equation}\label{eq:rays} 
\begin{aligned}
\(g_j(t),h_j(t)\) =&\left\{
\begin{aligned}
\( \frac{\sinh (\om_j t)}{\om_j},\cosh (\om_j t)\)\, ,\ &\text{ if
}\delta_j=-1\, ,\\ 
\(t,1\)\ ,\ &\text{ if }\delta_j=0 \, ,\\  
\(\frac{\sin (\om_j t)}{\om_j}, \cos(\om_j t)\)\, ,\ &\text{ if }
\delta_j=+1\, .
\end{aligned}
\right.
\end{aligned}
\end{equation}
Recall that if there exists $\delta_j=+1$, then $e^{-it{H}_V^\e}$ has
some singularities, periodically in time  
(see e.g. \cite{KR}). This affects the 
above formula with phase factors we did not write (which can be
incorporated in the definition of $(ig_j(t))^{1/2}$). 

Let us examine the Strichartz estimates we can deduce. If
$\delta_j=+1$ for any $j$, then only local in time Strichartz
estimates are available. This is not surprising, since the harmonic
oscillator has eigenvalues. On the other hand, if $\om_j=-1$ for any
$j$, then we have $\|U_V^\e(t)\|_{L^1\to L^\infty} \le |2\pi\e
t|^{-n/2}$, an estimate which is independent of all the $\om_j$'s: we
have the same dispersion as in the case with no potential (recall from
the previous section that the repulsive harmonic potential accelerates
the ``particle''), and global in time Strichartz estimates follow. Actually,
this property remains if $\delta_j=-1$ for \emph{at least one}
$j$. This stems from the formulation of Lemma~\ref{lem:strichartz}:
the function $\w(t)$ given by Mehler's formula is in the weak $L^1$
space for small values of $|t|$, and in the strong $L^1$ space for
``large'' values of $|t|$. 

\begin{rema}
We would like to point out that in the case where $\delta_j=-1$ for
at least one $j$, we have global in time Strichartz estimates, while
there are trapped trajectories. The one-dimensional case shows the
mechanism. If $n=1$ and $\om=1$, the classical symbol of $H_V^\e$ is
$\frac{1}{2}(\xi^2-x^2)$. Computing the Hamilton flow, we have
\begin{equation*}
x(t) = x_0 \cosh t +\xi_0 \sinh t = \frac{x_0+\xi_0}{2}e^{t} +
\frac{x_0-\xi_0}{2} e^{-t}\, .
\end{equation*}
If $x_0+\xi_0=0$, then the trajectory is trapped in the future. 
This aspect can be compared to the results of \cite{CKS95,Doi00},
where it is proved that smoothing effects occur in the future provided
that the classical trajectories are not trapped in the past. However,
the results of \cite{CKS95} include potentials which grow at most
linearly in $x$, and \cite{Doi00} does not consider the case of
potentials. On the other hand,
smoothing effects yield another approach to prove Strichartz estimates
(see e.g. \cite{ST02}). In our
case, there are trajectories trapped in the past or in the
future, but global in time Strichartz estimates are available. 
It seems
that the link between classical trajectories and (global in time)
Strichartz estimates remains to be clarified. 
\end{rema}

To introduce the tools we use in the rest of the analysis, recall the
approach followed in \cite{CaIHP}. We consider the case of an
isotropic harmonic potential, $V(x)=\frac{\om^2}{2}|x|^2$. In the
linear case $\lambda^\e=0$, $u^\e$ is therefore given by Mehler's
formula. A formal stationary argument (which can be justified) shows
that if $u_0^\e=f$ does not depend on $\e$, then for
$|t|<\frac{\pi}{2\om}$, 
\begin{equation}\label{eq:uapp}
u^\e(t,x)\Eq \e 0 \frac{1}{(\cos (\om t))^{n/2}}f\( \frac{x}{\cos (\om
t)}\) e^{-i\om \frac{|x|^2}{2\e}\tan (\om t)}=:u^\e_{\rm app}(t,x)\,
\, .
\end{equation}
In the nonlinear case, we expect that if, say $\lambda^\e =\O(\e^k)$
with $k$ sufficiently large, then nonlinear effects should be
negligible in the semi-classical limit, at least before the first
singularity at time $\frac{\pi}{2\om}$. This is proved in
\cite{CaIHP}, and we recall the argument in
Section~\ref{sec:semiclas}. A natural candidate for an approximate
solution is then given by $u^\e_{\rm app}$. To prove the approximation
is valid in the 
nonlinear case, $L^2$ estimates are not sufficient: $L^p$ estimates
are needed, for other values of $p$, and one can think of
Gagliardo--Nirenberg inequalities, $\|g\|_{L^p}\lesssim
\|g\|_{L^2}^{1-\delta(p)} \|\nabla g\|_{L^2}^{\delta(p)}$.   

However, these inequalities  yield 
$\|u^\e_{\rm app}(t)\|_{L^p} \lesssim \e^{-\delta(p)}$, which is
terrible in the limit $\e \to 0$ for $|t|<\frac{\pi}{2\om}$. This
``bad'' power of $\e$ stems from the highly oscillatory phase. Note
that the nonlinearity we consider, $F(z)=|z|^{2\si}z$, does not create
new harmonics in a single phase WKB asymptotics, if only one harmonic
is present initially. Therefore, in our case, there is only one phase
and one harmonic to take care about. This suggests to replace the
gradient in Gagliardo--Nirenberg inequalities by the operator
$\nabla_x\(e^{i\om \frac{|x|^2}{2\e}\tan (\om t)} \, \cdot \)$. In
view of the expression for $u^\e_{\rm app}$, we introduce more
precisely:
\begin{equation*}
J^\e(t) = i \cos (\om t) e^{-i\om \frac{|x|^2}{2\e}\tan (\om t)} 
\nabla_x\(e^{i\om \frac{|x|^2}{2\e}\tan (\om t)} \, \cdot \)= 
-\om^2 \sin(\om t)\frac{x}{\e}+i\cos(\om t)\nabla_x\, .
\end{equation*}
Introduce the operator 
\begin{equation*}
H^\e(t) =  
\cos(\om t)x+i\frac{\sin(\om t)}{\om}\e\nabla_x\, .
\end{equation*}
We can then rewrite the energy $E_V^\e$ (which is constant from
Lemma~\ref{lem:exlocV}), as
\begin{equation*}
E_V^\e = \frac{1}{2}\left\| \e J^\e(t)u^\e\right\|_{L^2}^2 +
\frac{\om^2}{2} \left\| H^\e(t)u^\e\right\|_{L^2}^2
+\frac{\lambda^\e}{\si +1} \left\| u^\e(t)\right\|_{L^{2\si +2}}^{2\si
+2}\, .
\end{equation*}
More generally, if the potential if of the form \eqref{eq:V}, we
define:
\begin{equation}\label{eq:JH}
\( 
\begin{array}{c}
J_j(t)\\
H_j(t)
\end{array}
\)= \( 
\begin{array}{cc}
-\delta_j\om_j^2 g_j(t) & h_j(t)/\e\\
\e h_j(t)& g_j(t)
\end{array}
\)
\( 
\begin{array}{c}
x_j/\e\\
i\e\d_j
\end{array}
\)\, ,\quad \forall j\ge 1\, .
\end{equation}
We set $J^\e =(J_j^\e)_{1\le j\le n}$ and $H^\e =(H_j^\e)_{1\le j\le
n}$. These operators have been known for years, in the linear 
theory (see e.g. \cite{Thirring,Robert}). They are \emph{Heisenberg
observables}:
\begin{equation*}
J^\e(t) = U_V^\e(t)i\nabla_x U_V^\e(-t)\quad ; \quad H^\e(t)=
U_V^\e(t)x U_V^\e(-t) \, .
\end{equation*}
The fact that such Heisenberg observables can be computed exactly is
due to the assumption that the potential is a second order
polynomial. Note that if $\delta_j=0$, then we recover the two
operators introduced in Section~\ref{sec:nls}: the gradient and the
Galilean operator. We will see in the next section that when $V$
is of the form \eqref{eq:V}, then $J^\e$ and $H^\e$ satisfy the same
properties as those emphasized in Section~\ref{sec:nls}. 
\section{Semi-classical analysis in a nonlinear framework}
\label{sec:semiclas}

The main features of the operators $J^\e$ and $H^\e$ are the
following:
\begin{lem}\label{lem:J&H}
Let $V$ be of the form \eqref{eq:V}. The operators $J^e$ and $H^\e$
satisfy:\\ 
\noindent $(i)$ They commute with the linear part of \eqref{eq:nlsp}:
\begin{equation*}
\left[ J^\e(t),i\e\d_t +\frac{1}{2}\e^2 \Delta -V(x)\right]=\left[
H^\e(t),i\e\d_t +\frac{1}{2}\e^2 \Delta -V(x)\right]= 0\, .
\end{equation*}
\noindent $(ii)$ They can be written as
\begin{align*}
J^\e_j(t) &=i h_j(t) e^{i\phi_1(t,x)/\e} \d_j \( e^{-i\phi_1(t,x)/\e}\,
\cdot\)\, ,\\ 
H^\e_j(t) &=i g_j(t) e^{i\phi_2(t,x)/\e} \d_j \(
e^{-i\phi_2(t,x)/\e}\, 
\cdot\)\, ,
\end{align*}
for some real-valued phases $\phi_1$ and $\phi_2$, solutions of the
eikonal equation \eqref{eq:eikonale}. \\
\noindent $(ii)'$ They yield weighted Gagliardo--Nirenberg
inequalities, such as
\begin{equation*}
\|f\|_{L^p}\lesssim \Big(\prod_{j=1}^n |h_j(t)|^{\delta(p)/n}\Big)
\|f\|_{L^2}^{1-\delta(p)}\|J(t)f\|_{L^2}^{\delta(p)}\, . 
\end{equation*}
\noindent $(ii)''$ If $F(z)=G(|z|^2)z$ is $C^1$, then $J(t)F(w)$ and
$H(t)F(w)$  are given by
\eqref{eq:Jder}. 
\end{lem}
The first point is trivial, since $J^\e$ and $H^e$ are
Heisenberg observables. The last two points follow from the second
one, which is the way we found the operator $J^\e$ in the case of the
isotropic harmonic potential. We can consider
that two algebraic miracles occur: first, we can compute explicitly
some interesting Heisenberg observables. Second, these operators can
be written like (ii), which is a nice property in view of
nonlinear problems. We discuss these aspects further into
details below. 

We now come to the issue of critical scales in the semi-classical
analysis of \eqref{eq:nlsp}. We describe two cases: 
\begin{itemize}
\item $V$ is an isotropic harmonic potential, and $u_0^\e=f$ does not
depend on $\e$.
\item $V$ is of the form \eqref{eq:V} and $u_0^\e$ is a concentrating
profile. 
\end{itemize}
The first case corresponds to the one that led us to introduce the
operators $J^\e$ and $H^\e$. In the linear case, due to the harmonic
potential, the solution focuses at the origin at time
$t=\frac{\pi}{2\om}$. More precisely, the geometry of the propagation
in the limit $\e\to 0$ is given by the Hamilton flow. The classical
Hamiltonian in this case is
$p(t,x,\tau,\xi)=\frac{1}{2}(|\xi^2|+\om^2|x|^2)$. The classical
trajectories (rays of geometric optics) are given by 
\begin{equation*}
x(t)=x_0 \cos (\om t)+\xi_0 \frac{\sin(\om t)}{\om}\, \cdot
\end{equation*}
Since $f$ does not depend on $\e$, there is no initial oscillation:
$\xi_0=0$, and rays meet at the origin, periodically in time. 
When $\lambda^\e=0$, one has the sharp estimate:
\begin{equation*}
\|u^\e(t)\|_{L^p} \lesssim \(\frac{1}{|\cos (\om t)|+\e |\sin(\om
t)|}\)^{\delta(p)}\, .
\end{equation*}
This follows easily from Lemma~\ref{lem:J&H} and the conservation of
the $L^2$ for solution to linear Schr\"odinger equations. Now assume
that $\lambda^\e = \e^{\alpha}$, for some $\alpha>0$. Then the
conservations of mass and energy (Lemma~\ref{lem:exlocV}) show that
the solution $u^\e$ is global in time. For $\alpha$ large, the linear
solution is expected to be a good approximation for the nonlinear
solution. More precisely, two r\'egimes must be considered. Before
focusing, the solution is of order $\O(1)$, and one can apply WKB
methods. This leads to a linear approximation if $\alpha>1$, nonlinear
otherwise. Next, at the focus, the linear solution is of order
$\e^{-n/2}$, more precisely, we infer from \eqref{eq:mehler}:
\begin{equation*}
u^\e_{\rm lin}\(\frac{\pi}{2\om},x \)\Eq \e 0 \(\frac{\om}{\e}\)^{n/2}
\widehat f\( \frac{\om x}{\e}\)\, , \quad \text{where }
\widehat f(\xi)=\int_{\R^n} e^{-ix.\xi}f(x)dx\, .
\end{equation*}
The nonlinear term can be viewed as a potential $\e^\alpha
|u^\e|^{2\si}$. Plugging the above asymptotics suggests that a
critical value for $\alpha$ at the focus is $\alpha =n\si$. 

With the arguments given in \cite{CaIHP}, it is possible to prove that
if $\alpha >\max(1,n\si)$, then the nonlinearity is negligible in the
limit $\e\to 0$, locally uniformly in time. The case $\alpha=n\si>1$
is shown to be critical:
\begin{theo}[\cite{CaIHP}]\label{theo:IHP}
Let $V(x)=\frac{|x|^2}{2}$, and $2 < r< \frac{2n}{n-2}$, with
$r=\infty$ if $n=1$. 
Assume that the nonlinearity
$z\mapsto |z|^{2\sigma}z $ is twice differentiable, and that
$\lambda^\e=\e^{n\si}$, with $n\si>1$. Assume moreover that 
$\si>\si_0(n)$.  
Then for $k\in \N$, the asymptotics of $u^\e$ for $\frac{\pi}{2} +
(k-1)\pi < t < \frac{\pi}{2} +k\pi$ is given, in $L^2\cap
L^r$, by:
$$u^\e(t,x) \Eq \e 0 
\frac{e^{in\frac{\pi}{4}-ink\frac{\pi}{2}}}{(2\pi |\cos
t|)^{n/2}}\widehat{S^k\psi_-}\left(\frac{-x}{\cos t} \right)
e^{-i\frac{|x|^2}{2\e}\tan t}\, ,$$
where $S^k$ denotes the $k$-th iterate of $S$ (which is well defined
under our assumptions on $f$), and $\psi_-(x):=
(2i\pi)^{-n/2}\widehat f (x)$. 
\end{theo}
Note that the phase shift $-n\frac{\pi}{2}$,
appearing at each focus crossing, is a linear phenomenon
(Maslov--Keller index, see \cite{Du}). 
Thus the only nonlinear effects at leading order occur at
the focuses, and are described, in average, by the scattering operator
$S$ associated to \eqref{eq:nls}. We give the main ideas of the
proof. The first step consists in proving that before the first focus
at time $t=\frac{\pi}{2}$, the function $u_{\rm app}^\e$ defined by
\eqref{eq:uapp} is a good approximation for the nonlinear solution
$u^\e$. This is equivalent to justifying a WKB asymptotics at leading
order, and relies on the sharp estimate given by the operator
$J^\e$ and the weighted Gagliardo--Nirenberg inequality it provides
(point (ii)' in Lemma~\ref{lem:J&H}). We can prove that $u_{\rm
app}^\e$ remains close to $u^\e$ (in a space which is essentially
$\Sigma$), up to time $\frac{\pi}{2}-\Lambda \e$, in the limit
$\Lambda \to +\infty$, that is, until focusing effects become relevant
at leading order. This shows that as predicted, the assumption $\alpha
>1$ makes the nonlinearity negligible outside the focus. 

When $t-\frac{\pi}{2}=\O(\e)$, $u^\e$
tends to concentrate near the origin, at scale $\e$; 
$V(x)u^\e$ becomes negligible, while the nonlinear potential
$\e^{n\si}|u^\e|^{2\si}$ is of order $\O(1)$. This can be proved
thanks essentially to the operator $H^\e$. Note that for
$t-\frac{\pi}{2}=\O(\e)$, the operator $J^\e$ and $H^\e$ can be
replaced  respectively by 
$\frac{x}{\e}+i\(t-\frac{\pi}{2}\)\nabla_x$ and  $\e \nabla_x$, 
which is another way to check that the potential is negligible. The scaling 
\begin{equation}\label{eq:scalingsemiclas}
u^\e(t,x)= \frac{1}{\e^{n/2}}\psi^\e \(
\frac{t-\frac{\pi}{2}}{\e}\virgp \frac{x}{\e} \)
\end{equation}
turns the description of the caustic crossing into a continuity issue
for a scattering problem. We have $\|U_0(-t)(\psi^\e(t)
-\psi(t))\|_{L^\infty(\R;\Sigma)} = o(1)$, where $\psi$ solves
\eqref{eq:nls}, with $\psi_-$ as a Cauchy datum at $t=-\infty$. This
scattering state stems from the transition between the two r\'egimes
for $u^\e$, which occurs for  $t= \frac{\pi}{2}-\Lambda \e$, in the
limit $\Lambda \to +\infty$. 

Next, $u^\e$ leaves the focus in a way described by the second point
of Proposition~\ref{prop:scattnls} and \eqref{eq:scalingsemiclas}. The
analysis is 
symmetric to the one before the focus; we prove the asymptotics until
$t=\pi$, where the situation is similar to that at time
$t=0$. Iterating the analysis yields Theorem~\ref{theo:IHP}.

\begin{rema}
With no additional effort, one can treat the case of initial plane
oscillations. Let $\xi_0\in \R^n$, and
introduce
$${\tt u}^\e(t,x)=u^\e(t,x-\xi_0 \sin t)e^{i\left( x-
\frac{\xi_0}{2}\sin t\right).\xi_0 \cos t /\e}\, .$$
If $u^\e$ solves \eqref{eq:nlsp} with $V(x)=\frac{|x|^2}{2}$ and
$u_0^\e =f$, then ${\tt u}^\e$ solves 
\eqref{eq:nlsp} with initial data ${\tt u^\e}_{\mid t=0}= f(x)
e^{i\frac{x\cdot \xi_0}{\e}}$. 
\end{rema}
The second case we analyze is when $V$ is of the general form
\eqref{eq:V} and  
\begin{equation}\label{eq:CICM}
u_0^\e(x)=\frac{1}{\e^{n/2}}\varphi\(\frac{x}{\e}\)\quad \(
\varphi\in \Sigma\)\quad ; \quad \lambda^\e =\e^{n\si}\, .
\end{equation}
If $\lambda^\e =\e^\alpha$ with $\alpha >\max
(1,n\si)$, then one can show that the evolution of $u^\e$ is linear,
at leading order. The case $\lambda^\e =-\e^{n\si}$ has been studied
by several authors (see e.g. \cite{BJ00,Keraani02}), 
in the case $\si<\frac{2}{n}$ and $\varphi(x)= R(x-x_0)e^{ix\cdot \xi_0}$,
where $R$ is a ground state. In that situation,
$u^\e$ evolves as a concentrating profile, with profile $R$,
along the Hamilton flow associated to $H_V$ with data
$(x_0,\xi_0)$. Heuristically, this is so because there is a 
balance between the dispersive effects associated to $H_V^\e$, and the
nonlinear effects (in the case $\lambda^\e<0$, the nonlinearity is
attractive). 

In the case $\lambda^\e = +\e^{n\si}$, the two effects mentioned above
tend to cumulate, and dispersion alters the shape of $u^\e$. Note that
even though $\lambda^\e>0$, global existence for $u^\e$ is not
obvious, since $V$ is not necessarily signed. 
\begin{theo}[\cite{CM}]\label{theo:Osaka}
Let $V$ of the form \eqref{eq:V}, and $u_0^\e$ given by
\eqref{eq:CICM}. Assume either that there exists 
$j$ such that $\delta_j \not = 1$, or 
that $\delta_j=1$ for all $j$ and the $\om_j$'s are not pairwise rationally
dependent. 
Suppose that the nonlinearity
$z\mapsto |z|^{2\sigma}z $ is twice differentiable. Then the following
holds. \\ 
1. For any $T>0$, there exists $\e(T)>0$ such that for $0<\e\leq
   \e(T)$, \eqref{eq:nlsp} has a unique solution $u^\e \in
   C([-T,T];\Sigma)$. \\
2. This solution satisfies the following asymptotics.
\begin{itemize}
\item For any $\Lambda >0$,
\begin{equation*}
\begin{aligned}
\limsup_{\e \to 0}&\sup_{|t|\leq \Lambda \e}\Big(
\left\|u^\e(t)-v^\e(t) \right\|_{L^2}+
\left\|\e\nabla_x u^\e(t)-\e\nabla_x v^\e(t) \right\|_{L^2}\\
& +
\left\|\frac{x}{\e}u^\e(t)- \frac{x}{\e}v^\e(t) \right\|_{L^2}\Big)
   =0\ ,\quad \text{where }v^\e(t,x)=\frac{1}{\e^{n/2}}\psi
   \left(\frac{t}{\e},\frac{x}{\e} \right),  
\end{aligned}
\end{equation*}
and $\psi\in C(\R;\Sigma)$ is the solution to \eqref{eq:nls} such that
$\psi_{\mid t=0}=\varphi$. 
\item Beyond this boundary layer, we have
\begin{equation*}
\begin{aligned}
\limsup_{\e \to 0}\sup_{\Lambda \e \leq \pm t\leq T}\Big(&
\left\|u^\e(t)-u_\pm^\e(t) \right\|_{L^2}+
\left\|\e\nabla_x u^\e(t)-\e\nabla_x u_\pm^\e(t) \right\|_{L^2}\\
&+\left\|x u^\e(t)-x u_\pm^\e(t) \right\|_{L^2}\Big)
   \Tend \Lambda {+\infty} 0,
\end{aligned}
\end{equation*}
where $u^\e_\pm \in C(\R;\Sigma)$ are the solutions to 
\begin{equation*}
i\e \d_t u_\pm^\e +\frac{1}{2}\e^2\Delta u_\pm ^\e =
V(x)u_\pm^\e\quad ;\quad 
u^\e_{\pm \mid t=0} = \frac{1}{\e^{n/2}}\psi_\pm
\left( 
\frac{x}{\e}\right),
\end{equation*}
and $\psi_\pm$ are given by Proposition~\ref{prop:scattnls}. 
\end{itemize}
\end{theo}
This result, as well as the analysis in \cite{CaIHP}, can be viewed as
a nonlinear analog to a result due to 
F.~Nier. In \cite{NierENS}, the author studies the
problem  
\begin{equation*}
i\e\d_t u^\e +\frac{1}{2}\e^2\Delta u^\e 
=V(x)u^\e +U\left( \frac{x}{\e}\right)u^\e\quad ;\quad
u^\e_{\mid t=0} = \frac{1}{\e^{n/2}}\varphi\left( \frac{x}{\e}\right),
\end{equation*}
where $U$ is a short range potential. The potential $V$ in that case
is bounded as well as
all its derivatives. Under
suitable assumptions, the influence of $U$ occurs near $t=0$ and is
localized near the origin, while only the value $V(0)$ of $V$ at the origin
is relevant in this r\'egime. For times $\e \ll |t|<T_*$, the
situation is different: the potential $U$ becomes negligible, while
$V$ dictates the propagation. Like in the nonlinear setting, the
transition between 
these two r\'egimes is measured by the scattering operator associated
to $U$. 

The assumption $\si\ge \si_0(n)$ 
 makes the nonlinear term short range. With
our scaling for the nonlinearity, this perturbation is relevant only
near the focus, where the potential is negligible, while
the opposite occurs for $\e\ll |t|\leq T$. 

The assumption on the potential $V$ in Theorem~\ref{theo:Osaka}  is
such that this result is 
complementary to the analysis in \cite{CaIHP}. It excludes the
phenomenon of ``total'' refocusing. Indeed, the rays of geometric
optics will all meet at one point for a time $t\not =0$ if and only if
$\delta_j=+1$ for all $j$ and the $\om_j$'s are pairwise rationally
dependent. As discussed in \cite{CaIHP}, this is the only case where
the nonlinear term could be relevant at leading order, past the
initial neighborhood of size $\e$. 
\begin{rema}
An explicit change of variable makes it possible to state the theorem
when $\varphi(\frac{x}{\e})$ is replaced by
$\varphi(\frac{x-x_0}{\e})e^{ix\cdot \xi_0/\e}$ 
(the phase factor $e^{ix\cdot \xi_0/\e}$ is actually not relevant, since
the only assumption we make is $\varphi\in\Sigma$), see \cite{CM}. 
\end{rema}

\begin{rema}
We did not mention the case $\alpha =1>n\si$. Due to the lack of
regularity of the nonlinearity $F(z)=|z|^{2\si}z$, this r\'egime can
be analyzed more easily when the power nonlinearity is replaced by a
Hartree-type nonlinearity. In \cite{CaMaSt}, we study the asymptotics
behavior as $\e\to 0$ of the solution $u^\e$ to 
\begin{equation*}
i\e \d_t  u^\e +\frac{1}{2}\e^2\Delta u^\e = \frac{|x|^2}{2}u^\e +
\e^\alpha \(|x|^{-\gamma}\ast |u^\e|^2\)u^\e \quad ; \quad u^\e_{\mid t=0} =
f\, ,
\end{equation*}
with $\gamma>0$, $\alpha \geq 1$ and $x\in\R^n$ ($n\geq 2$). In the
same spirit as above, the expected critical values for the parameters
are $\alpha>1$ or $\alpha=1$ to measure the relevance of the
nonlinearity outside of the focus, and $\alpha>\gamma$ or
$\alpha=\gamma$ to decide whether the nonlinearity plays an important
role at the caustic or not. When $\alpha =\gamma>1$,  results 
similar to Theorem~\ref{theo:IHP} are proved. When $\alpha=1>\gamma$,
then as 
expected, the nonlinearity is relevant only away from the focus (and
is measured by a slowly varying phase shift); the caustic crossing is
described by the Maslov--Keller index. When $\alpha=\gamma=1$, both
nonlinear effects mentioned above are expected, and the underlying
scattering operator would be a \emph{long range} scattering
operator. This case is treated only through a formal computation. 
\end{rema}

Theorems~\ref{theo:IHP} and \ref{theo:Osaka} show that when the
potential $V$ is of the form \eqref{eq:V}, the interactions between
the linear dynamics and the nonlinear effects can be understood quite
precisely. We mention two important limitations. First, we do not
treat super-critical cases, such as $\lambda^\e=\e^\alpha$ with
$\alpha <n\si$. The difficulty to treat such cases is not specific to
the presence of the potential $V$, however. Second, a natural
question is: what if $V$ is a general smooth sub-quadratic potential?
Some answers are given in \cite{CM}. As we already pointed out, the
compatibility between linear and nonlinear analysis in \cite{CaIHP,CM}
relies essentially in Lemma~\ref{lem:J&H}. In \cite{CM}, we ask: for
which potentials $V$ can an operator of the form
$\rho(t)e^{\phi/\e}\nabla_x ( e^{\phi/\e} \cdot)$ (compare with
Lemma~\ref{lem:J&H}, (ii)) commute with the linear part of
\eqref{eq:nlsp} (point (i) in Lemma~\ref{lem:J&H})? Easy computations
show that this is possible only if $V$ is a second order polynomial
(and when $\phi$ solves the eikonal equation \eqref{eq:eikonale}). For
initial data of the form \eqref{eq:CICM}, an interesting candidate as
a substitute to $J^\e$ (which saves the day away from the focus in
Theorems~\ref{theo:IHP} and \ref{theo:Osaka}) would 
be the Heisenberg observable 
$U_V^\e(t)\frac{x}{\e}U_V^\e(-t)$. Note that it coincides with $J^\e$ in the
case $V(x)=\frac{|x|^2}{2}$ if $\frac{\pi}{2}$ is taken as a new time
origin. Obviously, it satisfies the commutation property (i) in
Lemma~\ref{lem:J&H}. It also yields a weighted Gagliardo--Nirenberg
of the form stated in Lemma~\ref{lem:J&H}, (ii)', at least for $|t|\le
\delta$ small (the same as in
Proposition~\ref{prop:parametrice}). In Lemma~\ref{lem:J&H}, the
only interest of point (ii) is to imply (ii)' and (ii)'', so only the
action of this Heisenberg observable on nonlinearities of the form
$G(|z|^2)z$  remains to be understood. This issue seems to be
connected to Egorov theorem; we have no answer to provide.  

\section{Changing finite time blow-up}
\label{sec:blowup}

Guided by the semi-classical analysis described above, we can
understand the role of some potentials on the finite time blow-up
phenomenon. In the linear semi-classical
analysis, it is well known that the energy is carried by the
bicharacteristic curves. On the other hand, blow-up for \eqref{eq:nls}
is when the $L^2$-norm of $\nabla_x u$ becomes infinite, while the
$L^2$-norm of $u$ is constant; heuristically, the energy of $u$
concentrates ``somewhere''. One can expect that if in the linear case,
bicharacteristics meet, then the associated potential may encourage
finite time blow-up in the nonlinear (focusing) case; the two effects
cumulate. On the contrary, if the bicharacteristics spread out, then
the potential may compete with the attractivity of the
nonlinearity. We give several illustrations that provide a rigorous
support for these ideas. 

\begin{figure}[htbp]
\begin{center}
\input{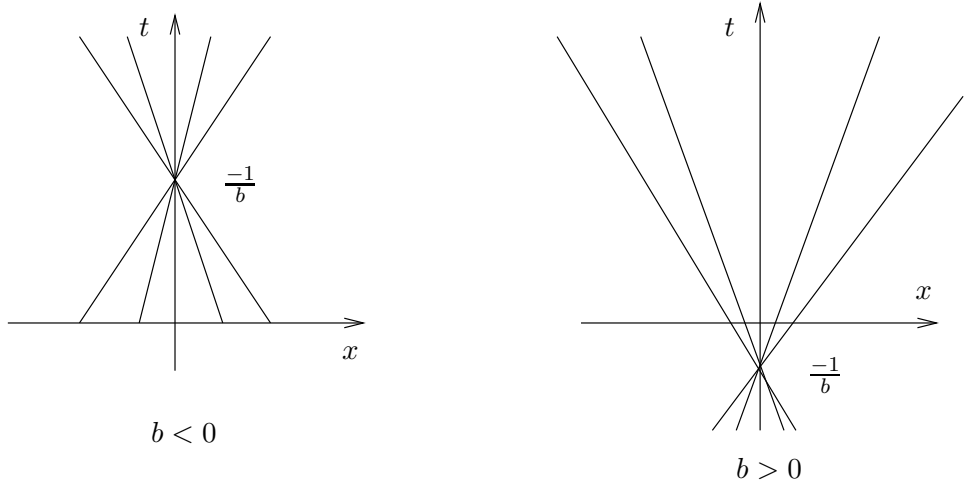}
\caption{Geometry of rays: case of quadratic oscillations.}
\label{fig:quad}
\end{center}
\end{figure}

The first examples follow from the same heuristic, when one plays with
the phase of the initial data instead of a potential in the
equation. Consider the linear equation
\begin{equation*}
i\e\d_t u^\e +\frac{1}{2}\e^2\Delta u^\e = 0\quad ;\quad u^\e(0,x) =
f(x)e^{ib\frac{|x|^2}{2\e}}\, ,\ b\in \R\, .
\end{equation*}
The Hamilton flow is given by
\begin{equation*}
\dot t =1\quad ;\quad \dot x=\xi\quad ;\quad \dot \xi =\dot\tau= 0\quad ;\quad
 x_{\mid t=0}=x_0\quad ;\quad
\xi_{\mid t=0} = b x_0\, ,
\end{equation*}
and the classical trajectories are: $x(t) = x_0(1+bt)$. In particular
$x(\frac{-1}{b}) =0$ for all $x_0\in \R$; if $b<0$, then rays meet at
the origin in the future, while if $b>0$, they met at the origin in
the past (see Figure~\ref{fig:quad}). 
Following the above heuristic discussion, one may expect
that quadratic initial oscillations alter the blow-up phenomenon in the
nonlinear case. 
This is proved in \cite{CW92}: let $u$ solve \eqref{eq:nls} with
$u_0\in \Sigma$, and suppose that blow-up may occur, that is
$\lambda<0$ and $\si\ge \frac{2}{n}$. Denote $u^b$ the solution of
\eqref{eq:nls} with initial datum
$u_0^b(x)=u_0(x)e^{ib\frac{|x|^2}{2}}$. Then for $b\gg 1$, $u^b$ is
global \emph{in the future} ($t\ge 0$). On the other hand, if
$E_0(u_0)<0$, then we know from 
Proposition~\ref{prop:glassey} that $u$ blows up at time $T>0$, say. It
is proved in \cite{CW92} that if $b<0$, then $u_b$ blows up at 
time $T^b\le \frac{-1}{b}$. For $b<\frac{-1}{T}$, the blow up phenomenon
occurs sooner for $u^b$ than for $u$. In the conformal case
$\si=\frac{2}{n}$, the critical values for $b$ can be explicitly
related to the blow-up time $T$ (see \cite{CW92,CazCourant}).

Consider now the  equation
\begin{equation*}
i\e\d_t u^\e +\frac{1}{2}\e^2\Delta u^\e = (E\cdot x)
 u^\e\quad ;\quad u^\e(0,x) =
f(x)e^{ib\frac{|x|^2}{2\e}}\, ,\ b\in \R\, .
\end{equation*}
The Hamilton flow is given by
\begin{equation*}
\dot t =1\quad ;\quad \dot x=\xi\quad ;\quad \dot \xi = -E
\quad ;\quad \dot\tau= 0\quad ;\quad
 x_{\mid t=0}=x_0\quad ;\quad
\xi_{\mid t=0} = b x_0\, ,
\end{equation*}
and the classical trajectories are: $x(t) =
x_0(1+bt)-\frac{t^2}{2}E$. Rays meet at $\underline{x} =
-\frac{1}{2b^2}E$ at time $\frac{-1}{b}$. This is the same phenomenon
as above, shifted in space. 
Recall that we saw in Section~\ref{sec:stark} that the
introduction of a linear potential $E\cdot x$ does not change the time
of blow-up, but only shifts the set where this occurs. Here again,
intuition and results meet. 

\begin{figure}[htbp]
\begin{center}
\input{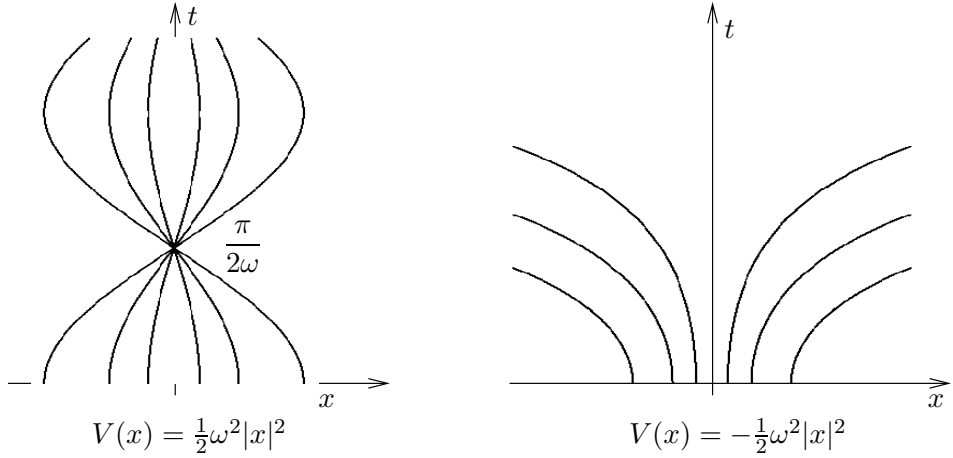}
\caption{Geometry of rays: isotropic quadratic potential.}
\label{fig:curv}
\end{center}
\end{figure}

The last two cases we consider are isotropic harmonic potential and
isotropic repulsive harmonic potential, with no initial rapid
oscillation:
\begin{equation*}
i\e\d_t u^\e_\pm +\frac{1}{2}\e^2\Delta u^\e_\pm = \pm\om^2\frac{|x|^2}{2}
 u^\e_\pm\quad ;\quad u^\e_\pm(0,x) =
f(x)\, .
\end{equation*}
We have $x(t) = x_0 h(t)$, where $h$ is given by
\eqref{eq:rays}. Thus, $x_+(t)= x_0 \cos (\om t)$, and $x_-(t)=
x_0\cosh (\om t)$. This is illustrated in Figure~\ref{fig:curv}. Note
the analogy with Figure~\ref{fig:quad}.  
We already saw in Section~\ref{sec:isoharmocrit} that the introduction
of an isotropic harmonic potential may anticipate the blow-up time
when $\si=\frac{2}{n}$,
just like quadratic oscillations in \cite{CW92}. Similarly, the
isotropic repulsive harmonic potential delays, or prevents, finite time
blow-up. 

To complete the picture, we have to study the case $\lambda<0$,
$\si>\frac{2}{n}$. 
\begin{theo}[\cite{CaAHP,CaSIMA}]\label{th:blow}
 Let $u_0 \in \Sigma$, $\lambda <0$, $\si \ge \frac{2}{n}$ and $\si
 <\frac{2}{n-2}$ 
 if $n\geq 3$. Let $u_\pm$ be the solutions of \eqref{eq:nlsp}
 with $\e=1$ and $V_\pm(x)=\pm \om^2\frac{|x|^2}{2}$.  \\
$1.$ If $E_{V_+}(u_0) \le \frac{1}{2}\om^2\|x u_0\|_{L^2}^2$, then $u_+$
 blows up at time $T_+^\om\le \frac{\pi}{2\om}$. \\
$2.$ If the initial datum $u_0$ satisfies
$$ \frac{1}{2}\|\nabla u_0\|_{L^2}^2
+\frac{\lambda}{\si +1}\| 
u_0\|_{L^{2\si +2}}^{2\si +2}<-\frac{\om^2}{2}\|x u_0\|_{L^2}^2,$$
then $u_-$ blows up in finite time, in the
future or in the past.\\
$3.$ If the initial datum $u_0$ satisfies
$$ \frac{1}{2}\|\nabla u_0\|_{L^2}^2
+\frac{\lambda}{\si +1}\| 
u_0\|_{L^{2\si +2}}^{2\si +2}<-\frac{\om^2}{2}\|x u_0\|_{L^2}^2-\om 
\left|\operatorname{Im}\int \overline{u_0}x\cdot \nabla_x
u_0 \right|,$$
then $u_-$ blows up in finite time, in the
future \emph{and} in the past.\\
$4.$ There exists $\omega_1>0$ such that for any $\om \geq  \omega_1$,
the solution $u_-$ is global in time. 
\end{theo}
\begin{rema}
Notice that in the first point, finite time blow-up occurs for a range
of positive values of the energy (it was known that if $E_{V_+}<0$,
then finite time blow-up occurs, see \cite{CazCourant}). This is in
sharp contrast with the 
case of \eqref{eq:nls}, where zero energy solutions may be global. 
When $\si =\frac{2}{n}$ the solitary wave $e^{i t}R(x)$, where $R$ is
the ground  state, solves \eqref{eq:nls}, is global in time and has
zero energy.  Note that the condition in the first point also reads
$E_0(u_0)\le 0$ (but $E_0$ is not the energy associated to that
equation!). 
\end{rema}
The proof of points 2 and 3 relies on the Zakharov--Glassey method,
just like Proposition~\ref{prop:glassey}, so we shall not discuss it,
and refer to \cite{CaSIMA}. The proofs of points 1 and 4 rely on two
conservation laws, which are more precise than the conservation of
energy, and can be viewed as analogs to the pseudo-conformal
conservation law \eqref{eq:pseudoconf}. These laws have a geometric
meaning, since they ``follow'' the propagation.  
\begin{lem}\label{lem:pseudoharmo}
Let $u_0 \in \Sigma$, and $\si
 <\frac{2}{n-2}$ 
 if $n\geq 3$. Let $u_\pm$ be the solutions of \eqref{eq:nlsp}
 with $\e=1$ and $V_\pm(x)=\pm \om^2\frac{|x|^2}{2}$.  Introduce
\begin{align*}
    E_+^1(t)&:=\frac{1}{2}\|J_+(t) u_+\|_{L^2}^2+\frac{\lambda}{\si
      +1}\cos^2(\om t)
\|u_+(t)\|_{L^{2\si +2}}^{2\si +2}\,  ,\\
E_+^2(t)&:=\frac{\om^2}{2}\| H_+(t) u_+\|_{L^2}^2+\frac{\lambda}{\si
      +1}\sin^2(\om t)
\|u_+(t)\|_{L^{2\si +2}}^{2\si +2}\, ,\\  
E_-^1(t) &:= \frac{1}{2}\|J_-(t)u_-\|_{L^2}^2 +\frac{\lambda}{\si
+1}\cosh^2(\om t)\|u_-(t)\|_{L^{2\si+2}}^{2\si+2}\, , \\
E_-^2(t) &:= \frac{-\om^2}{2}\|H_-(t)u_-\|_{L^2}^2 -\frac{\lambda}{\si
+1}\sinh^2(\om t)\|u_-(t)\|_{L^{2\si+2}}^{2\si+2}\, , 
\end{align*}
where $J_+$ stands for $J$ when $V(x)=+\om^2\frac{|x|^2}{2}$, and so
on. Note that $E_\pm^1(t)+ E_\pm^2(t)\equiv E_{V_\pm}$. 
We have:
\begin{align}
    \frac{d E_+^1}{dt}=  \frac{\om \lambda}{2\si
    +2}(n\si-2) \sin(2\om t)\|u_+(t)\|_{L^{2\si +2}}^{2\si +2}\,
    ,\label{eq:E1+}\\
    \frac{d E_-^1}{dt}= \frac{\om \lambda}{2\si
    +2}(2-n\si) \sinh (2\om t)\|u_-(t)\|_{L^{2\si +2}}^{2\si +2}\,
    .\label{eq:E1-} 
\end{align}
\end{lem}
\begin{rema}
These algebraic results can be proved in a classical way. Notice that
one turn $V_+$ into $V_-$ by replacing $\om $ with $i\om$ (and {\it vice
versa}). 
\end{rema}
The first point of Theorem~\ref{th:blow} follows easily. 
Assume that
$E_{V_+}(u_0)\le \frac{\om^2}{2}\|x u_0\|_{L^2}^2$, and suppose that
$u_+$ exists in $\Sigma$ 
up to time $\frac{\pi}{2\om}$. Then $E_+^1(0)=E_0(u_0)\le 0$, and from
\eqref{eq:E1+}, if $\si \ge \frac{2}{n}$ and $\lambda <0$, 
\begin{equation*}
\frac{d E_+^1}{dt} \leq 0\ , \ \ \forall t \in \left[0,\frac{\pi}{2\om}
\right]\, . 
\end{equation*}
This implies
\begin{equation*}
E_+^1\left(\frac{\pi}{2\om} \right) \le 0\, . 
\end{equation*}
But from the definition of $E_+^1$,  
\begin{equation*}
E_+^1\left(\frac{\pi}{2\om} \right) =  \frac{\om^2}{2}\left\| x u_+
\left(\frac{\pi}{2\om} \right) \right\|^2_{L^2}\, ,
\end{equation*}
so this leads to a contradiction (unless $u_+\equiv 0$, which
means that $u_0\equiv 0$). Therefore, $u_+$ does not remains in
$\Sigma$ up to time $\frac{\pi}{2\om}$. It is easy to conclude that this is
so because there exists $T\le\frac{\pi}{2\om}$ such that
\begin{equation*}
\lim_{t\to T} \left\|\nabla_x u_+ (t)\right\|_{L^2} =\infty\, .
\end{equation*}
Notice that the larger $\om$, the sooner the blow-up. Note that the sufficient
condition to have finite time blow-up does not depend on the value of
$\om>0$ though. 

Even if the conditions 2 and 3 in Theorem~\ref{th:blow} become void
as $\om\to +\infty$, this does not mean that the last point is
true. In the conformal case $\si=\frac{2}{n}$, the last point is
explicit, as we saw in Section~\ref{sec:isorepharmocrit}, up to a
characterization of finite time blow-up in this case where the
potential is unbounded from below. Roughly speaking, the
characterization of global existence is the same as in
Theorem~\ref{theo:nlslocal}. The energy $E_{V_-}$ is the sum
of three terms, and the term corresponding to the nonlinearity is
controlled by the gradient and the $L^2$ norm of the solution, which
is conserved. Thus, if the gradient remains bounded in $L^2$, then each
term of the energy remains bounded. Now the determinant of the matrix
in \eqref{eq:JH} 
is constant, equal to $-1$. Finally, to prove global
existence, it is sufficient to check that the $L^2$ norm of $J_-(t)u$
remains bounded.  

We briefly sketch the proof of the last point of Theorem~\ref{th:blow},
and refer to \cite{CaSIMA} for details. The first step 
consists in noticing that the usual method to prove local existence
still works for \eqref{eq:nlsp} with $\e=1$ and
$V(x)=-\frac{\om^2}{2}|x|^2$, and does not ``see'' the parameter
$\om\ge 0$. Indeed, as we already noticed, Mehler's formula yields the
same dispersive estimate as in the case $\om=0$, hence the same
Strichartz estimates. Moreover, the
operator $J_-$ plays a role analog to that of $\nabla_x$ in the case
of \eqref{eq:nls}: it commutes with the linear part of the equation, acts
on the nonlinearity like a derivative, and yields weighted
Gagliardo--Nirenberg estimates with a weight uniformly bounded in
$\om\ge 0$ (see Lemma~\ref{lem:J&H}, (ii)', and recall that $\cosh
x\ge 1$). Therefore, there exists $t_0>0$ independent of $\om\ge 0$
such that $u_-\in C(]-2t_0,2t_0[;\Sigma)$ and $A(t)u$ is bounded in
$L^\infty([-t_0,t_0]:L^2)$ uniformly in $\om\ge 0$ for any $A\in
\{Id,J_-,H_-\}$. Now integrate \eqref{eq:E1-} between $t_0$ and
$t>t_0$. 
Since $\lambda (2-n\si)\ge 0$, 
\begin{equation*}
\begin{aligned}
E_-^1(t)&\le E_-^1(t_0) + C\om 
\int_{t_0}^t \sinh (2\om s)\|u_-(s
)\|_{L^{2\si +2}}^{2\si +2} ds\\
&\le C(C_0) + C\om \int_{t_0}^t \frac{\sinh (2\om s)}{\cosh(\om
s)^{n\si}} \|J_-(s)u_-\|_{L^2}^{n\si}ds\, . 
\end{aligned}
\end{equation*}
The constants in the last estimate do not depend on
$\om$. Define
\begin{equation*}
y(t):= 
\dis \sup_{t_0\leq s\leq t}\|J_-(s)u_-\|_{L^2}^2\,  . 
\end{equation*}
We have 
\begin{equation*}
E_-^1(t)\leq C(C_0)+ Cy(t)^{n\si/2}\int_{t_0}^t \frac{\sinh (2\om
s)}{\cosh(\om s)^{n\si}}ds
\leq C(C_0) + \frac{C}{\cosh (\om t_0)^{n\si -2}}y(t)^{n\si/2}\ .
\end{equation*}
We finally have
\begin{equation*}
y(t) \leq C(C_0) +\frac{C}{\cosh(\om t_0)^{n\si
-2}}y(t)^{n\si/2}\ .
\end{equation*}
For $\si>\frac{2}{n}$, we conclude by a bootstrap argument, since the constants
and $t_0$ do not depend on $\om>0$, and
\begin{equation*}
\cosh(\om t_0)\Tend \om {+\infty} +\infty \, .
\end{equation*}
This yields a uniform bound for $\|J_-(t)u_-\|_{L^2}$, and proves
global existence. 

This proof relies on the evolution law \eqref{eq:E1-}, which seems to
be bound to the case of \emph{isotropic} potentials. A possible
question is to ask whether a similar result holds when $V$ is
of the form \eqref{eq:V}, with, say, $\delta_1=-1$ (recall that the
$V$ is non-negative, things are rather well understood). An answer is given
in \cite{CaQuad}:
\begin{theo}[\cite{CaQuad}]\label{theo:quad}
Take $\e=1$. Let $n\ge 2$, $\lambda\in \R$, $\si\ge \frac{2}{n}$ with
$\si<\frac{2}{n-2}$ if $n\ge  
3$, and $u_0\in \Sigma$. Let $V$ be of the form \eqref{eq:V} with
$\delta_1=-1$, and 
denote
\begin{equation*}
\om_\pm = \max \left\{ \om_j\ ;\ \delta_j=\pm 1\right\}\quad (\om_+
=0\text{ if there is no }\delta_j=+1)\, .
\end{equation*}
Then there exists $\Lambda=\Lambda (n,
\si, |\lambda|,\|u_0\|_\Sigma)$ such that for 
\begin{equation}\label{eq:condquad}
\om_- \ge \Lambda (1+\om_+)+ \frac{2\si^2}{2-(n-2)\si}(1+\om_+)\ln
(1+\om_+)\, , 
\end{equation}
the solution $u$ to \eqref{eq:nlsp}  is global in time, $u\in
C(\R;\Sigma)$.
\end{theo}
The statement can be summarized as follows: if the repulsive
force is sufficiently strong compared to other effects (linear
confinement is overcome if $\om_-\gg \om_+$, nonlinear effects are 
overcome if $\om_-\gg 1$), then the solution is global.
The strategy of the proof is as follows. First, in the
same spirit as in \cite{CaSIMA}, we analyze
the local existence result, to bound from below the local existence
time, in term of the parameters $\om_j$. Then, we notice
that we obtain a time 
at which $u$ is defined and small, and for which therefore
the nonlinearity is not too strong. We consider the solution of
the linear equation (\eqref{eq:nlsp} with $\lambda=0$) that coincides
with $u$ at that time. A continuity argument shows that $u$ cannot move
away too much from this linear solution. Since the linear solution is
global, so is $u$. Note that the nonlinear term in $\om_+$ in
\eqref{eq:condquad} is zero when $\om_+=0$, that is when there is no
confinement. 

\section{More on global existence}
\label{sec:global}

In the last point of Theorem~\ref{th:blow}, and in
Theorem~\ref{theo:quad}, we saw that if we consider a quadratic
potential which is sufficiently repulsive, then the solution to
\eqref{eq:nlsp} is global. In that case, we even have scattering: the
solution to \eqref{eq:nlsp} is asymptotically linear, as time becomes
infinite. 

It is proved in \cite{CaSIMA} that is the nonlinearity is defocusing
($\lambda>0$), and the potential $V$ is the isotropic repulsive
harmonic potential, then the solution $u$ of \eqref{eq:nlsp} is global, $u\in
C(\R;\Sigma)$. This follows from \eqref{eq:E1-}. We then have:
\begin{prop}[\cite{CaSIMA,CaQuad}]\label{prop:scatt}
 Let $\e=1$, $\lambda,\si >0$, with $\si <\frac{2}{n-2}$
 if $n\geq 3$.   \\ 
$1.$ Assume $V(x)=-\frac{\om^2}{2}|x|^2$.  
\begin{itemize}
\item For every $u_- \in \Sigma$, there exists a unique $u_0\in \Sigma$
 such that the maximal solution $u\in C(\R;\Sigma)$ to
 \eqref{eq:nlsp} satisfies
$$\left\| U_V(-t)u(t)-u_-\right\|_\Sigma \Tend
t {-\infty} 0\, .$$
\item For every $u_0\in \Sigma$, there exists a unique $u_+\in \Sigma$ 
 such that the maximal solution $u\in C(\R;\Sigma)$ to
 \eqref{eq:nlsp} satisfies
$$\left\| U_V(-t)u(t)-u_+\right\|_\Sigma \Tend
t {+\infty} 0\, .$$
\end{itemize}
$2.$ Suppose that $V$ is of the form \eqref{eq:V} with
$\delta_1=-1$.
\begin{itemize}
\item For every $u_- \in \Sigma$, there exist $T$ finite and a unique
$u\in C(]-\infty,T];\Sigma)\cap L^q\( ]-\infty,T] ;
L^{2\si+2}(\R^n)\)$,
where $q=\frac{4\si+4}{n\si}$,   solution to
 \eqref{eq:nlsp} such that
$$\left\| U_V(-t)u(t)-u_-\right\|_\Sigma \Tend
t {-\infty} 0\, .$$
\item Let $u_0\in \Sigma$. Then taking $\Lambda$ larger
in \eqref{eq:condquad} if necessary, there exists a unique $u_+\in \Sigma$ 
 such that the maximal solution $u\in C(\R;\Sigma)$ to
 \eqref{eq:nlsp} satisfies
$$\left\| U_V(-t)u(t)-u_+\right\|_\Sigma \Tend
t {+\infty} 0\, .$$
\end{itemize}
\end{prop}
Notice that unlike in Proposition~\ref{prop:scattnls}, there is no
additional assumption on $\sigma$, simply $\si>0$: all power-like
nonlinearities are short range for the Hamiltonian $H_V$ under our
hypotheses. 
 
The key ingredient of the proof consists in noticing that in the two
cases, global existence follows from the boundedness of
$J(t)u$ in $L^\infty(\R;L^2)$. The weighted Gagliardo--Nirenberg
inequality of
Lemma~\ref{lem:J&H} yields exponential decay for $u$, and any positive
power of
an exponentially decreasing functions is integrable at infinity. 

The restriction $u\in C(]-\infty,T];\Sigma)$ in the first part of
point 2  comes from the fact that in
general, we cannot prove that even 
if $\lambda>0$, the solution $u$ is defined globally in time. 

\bibliographystyle{amsplain}
\bibliography{../../Caustics/carles}

\end{document}